\documentclass[reqno,12pt]{article}
\usepackage{a4wide}

\usepackage{amsmath}
\usepackage{amssymb}

\newtheorem{thm}{Theorem}[section]
\newtheorem{lemma}[thm]{Lemma}

\def\tanh{{\rm tanh}}

\def\P{{\mathbb P}}
\def\R{{\mathbb R}}
\def\Z{{\mathbb Z}}
\def\N{{\mathbb N}}

\def\W{{\cal W}}

\def\C{{\C}}
\def\C{{\cal D}}

\def\sqr{\vcenter{
         \hrule height.1mm
         \hbox{\vrule width.1mm height2.2mm\kern2.18mm\vrule width.1mm}
         \hrule height.1mm}}                  
\def\square{\ifmmode\sqr\else{$\sqr$}\fi}

\def\la{{\lambda}}

\def\a{{\alpha}}

\def\proof{\noindent{\bf Proof. }}

\def\C{{\bf C}}

\def\be{\begin{equation}}
\def\ee{\end{equation}}

\def\Th{\vskip 2mm\noindent \bf Theorem \,}
\def\Re{\vskip 2mm\noindent{\bf Remark}:\,\it }
\def\Dp{\displaystyle}
\def\ba{\begin{array}}
\def\ea{\end{array}}

\usepackage{hyperref}

\begin{document}

\title{\sf Two-dimensional Poisson Trees converge to the Brownian web}
\author{P. A. Ferrari$^{~\diamondsuit}{}^*$
   L. R. G.
Fontes$^{~\diamondsuit}$\thanks{\scriptsize Partially supported by CNPq grants
  520811/96-8, 300576/92-7
and 662177/96-7 (PRONEX) and FAPESP grant 99/11962-9}\ \ \ \ \ Xian-Yuan
Wu$^{~\heartsuit}$\thanks{\scriptsize Research supported in part by Foundation
of Beijing Education Bureau and FAPESP through grant 01/02577--6.
} }
\date{}
\maketitle

\begin{center}\vskip-1cm
\begin{minipage}{13cm}
{\scriptsize\noindent $^\diamondsuit$ Instituto
de Matem\'atica e Estat\'{\i}stica, Universidade de S\~ao Paulo\\
\mbox{}\,\,\,\,\, Rua do Mat\~ao 1010, 05508-090 S\~ao Paulo
- SP, Brasil \vskip 1mm \noindent\hskip 3.7mm{E-mails:
\texttt{pablo@ime.usp.br}, \texttt{lrenato@ime.usp.br}} \vskip 4mm
\noindent{$^\heartsuit$ Department of Mathematics, Capital Normal
University, Beijing, 100037, China }

\noindent \hskip 3.7mm 
{E-mail: \texttt{wuxy@mail.cnu.edu.cn}, \texttt{xianyuan\_wu@hotmail.com}}}
\end{minipage}
\end{center}

\section*{}
\paragraph{Abstract}
The \emph{Brownian web} can be roughly described as a family of coalescing
one-dimensional Brownian motions starting at all times in $\R$ and at all 
points of $\R$. It was introduced by Arratia; a variant was then studied by 
T\'oth and Werner; 
another variant was analyzed recently by Fontes, Isopi, Newman and Ravishankar.
The two-dimensional \emph{Poisson tree} is a family of continuous 
time one-dimensional random walks with uniform jumps in a bounded interval.  
The walks start at the space-time points of a homogeneous Poisson process in
$\R^2$ and are in fact constructed as a function of the point process. This
tree was introduced by Ferrari, Landim and Thorisson. By verifying criteria 
derived by Fontes, Isopi, Newman and Ravishankar, we show that, when properly 
rescaled, and under the topology introduced by those authors, Poisson trees 
converge weakly to the Brownian web.

\paragraph{Keywords and phrases} 
Brownian web, Poisson tree, Donsker's invariance principle, coalescing
one-dimensional Brownian motions

\paragraph{1991 Mathematics Subject Classification} 60K35, 60F17


\section{Introduction and results}
\renewcommand{\theequation}{1.\arabic{equation}}
\setcounter{equation}{0}
\paragraph{The Two-dimensional Poisson Tree}Let $S$ be a two-dimensional
homogeneous Poisson process of parameter $\lambda$. $S$ is a random subset of
$\R^2$, $s\in S$ has coordinates $s_1,s_2$.

For $x=(x_1,x_2)\in \R^2$, $t\ge x_2$ and $r>0$, let $M(x,t,r)$ be
the following rectangle 
\begin{equation}
  \label{1'}
  M(x,t,r):=
\{(x'_1,x'_2)\,:\, \, |x'_1-x_1|\le r,\;x_2\le x'_2\le t\}. 
\end{equation}
As
$t$ grows, the rectangle gets longer. The first time $t$ that
$M(x,t,r)$ hits some (or {\it another}, when $x\in S$) point of
$S$ is called $\tau(x,S,r)$; this is defined by \begin{equation}\label{2'}
\tau(x,S,r):= \inf\{t>x_2\,:\,M(x,t,r)\cap
(S\setminus\{x\})\neq\phi\}. \end{equation} The hitting point is the point
$\alpha (x)\in S$ defined by \begin{equation}
  \label{1}
  \alpha (x)\;:=\; M(x,\tau(x,S,r),r)\cap (S\setminus\{x\}),
\end{equation}
which consists of a unique point almost surely. If $x=$ some
$s\in S$, we say that $\a(x)=\a(s)$ is the \emph{mother} of $s$
and that $s$ is a \emph{daughter} of $\alpha(s)$. Let $\alpha^0(x)
= x$ and iteratively, for $n\ge 1$, $\alpha^{n}(x) =
\alpha(\alpha^{n-1}(x))$. For the case of $x=$ some $s\in S$,
$\a^n(x)=\a^n(s)$ is the $n$th grand mother of $s$.

Now let $G=(V,E)$ be the random directed graph with vertices $V=S$ and edges
$E=\{(s,\a(s))\,:\, s\in S\}$. This model was proposed by Ferrari, Landim and
Thorisson \cite{2} who also proved that $G$ has a unique connected component. It
is straightforward to see that there are no circuits in the graph $G$ which is
then called the two-dimensional \emph{Poisson tree}.

This model is related to a model of drainage networks of Gangopadhyay, Roy 
and Sarkar \cite{grs}, which can be viewed as a discrete space, long range 
version of the two-dimensional Poisson tree.

\paragraph{Brownian web} 
Arratia \cite{3,4}, and later T\'oth and Werner \cite{5} constructed random
processes that formally correspond to coalescing one-dimensional Brownian
motions starting from every space-time point. Fontes, Isopi, Newman and
Ravishankar \cite{1} (see also \cite{1a}) extended their work by constructing
and characterizing the so-called \emph{Brownian web} as a random variable
taking values in a metric space whose points are compact sets of paths. The
space is given as follows.

Let $(\bar \R^2,\rho)$ be the completion (or compactification) of
$\R^2$ under the metric $\rho$ defined as 
\begin{equation}\label{2}
\rho(x,y)=\left
|\frac{\tanh(x_1)}{1+|x_2|}-\frac{\tanh(y_1)}{1+|y_2|}\right |\vee
|\tanh(x_2)-\tanh(y_2)|, 
\end{equation} 
for any $x=(x_1,x_2), y=(y_1,y_2)\in
\bar\R^2: 
=[-\infty,\infty]^2$, namely, the image of
$\bar\R^2$ under the mapping
\begin{equation}\label{3}
x=(x_1,x_2)\longmapsto
(\Phi(x_1,x_2),\Psi(x_2))\equiv\left 
(\frac{\tanh(x_1)}{1+|x_2|},\tanh(x_2)\right ).
\end{equation}

For $t_0\in [-\infty,\infty]$, let $C[t_0]$ denote the set of
functions $f$ from $[t_0,\infty]$ to $[-\infty,\infty]$ such that
$\Phi(f(t),t)$ is continuous. Define
\begin{equation}\label{4}\Pi=\bigcup_{t_0\in[-\infty,\infty]}C[t_0]\times
\{t_0\},\end{equation} where $(f,t_0)\in \Pi$ then represents a path in $\bar
\R^2$ starting at $(f(t_0),t_0)$.

For $(f(t_0),t_0)\in\Pi$, we denote by $\hat f$ the function that
extend $f$ to all $[-\infty,\infty]$ by setting it equal to
$f(t_0)$ for all $t<t_0$. We take
\begin{equation}\label{5}d((f_1,t_1),(f_2,t_2))=(\sup_t|\Phi(\hat
f_1(t),t)-\Phi(\hat f_2(t),t)|)\vee|\Psi(t_1)-\Psi(t_2)|,\end{equation} then
we get a metric space $(\Pi,d)$ of paths with specified starting
points in space-time. It is straightforward to check that
$(\Pi,d)$ is a complete separable metric space.

Let now $\cal H$ denote the set of compact subset of $(\Pi,d)$,
with $d_{\cal H}$ the induced Hausdorff metric, i.e., \begin{equation}
\label{6} d_{\cal H}(K_1,K_2)=\sup_{g_1\in K_1}\inf_{g_2\in
K_2}d(g_1,g_2)\vee\sup_{g_2\in K_2}\inf_{g_1\in K_1}d(g_1,g_2).
\end{equation} $({\cal H},d_{\cal H})$ is also a complete separable metric
space. Denote by ${\cal F}_{\cal H}$ the corresponding Borel
$\sigma$-algebra generated by $d_{\cal H}$.

The Brownian web is characterized in \cite{1} as a $({\cal
H},{\cal F}_{\cal H})$-valued random variable $\bar \W$ (or its
distribution $\mu_{\bar \W}$). Define the \emph{finite-dimensional
distributions} of $\bar \W$ as the induced probability measures
$\mu_{(x_1,t_1;\ldots;x_n,t_n)}$ on the subsets of paths starting
from any finite deterministic set of points
$(x_1,t_1),\ldots,(x_n,t_n)$ in $\R^2$.

Given $t_0\in\R$, $t>0$, $a<b$, and a $({\cal H},{\cal F}_{\cal
H})$-valued random variable $V$, let $\eta_{_V}(t_0,t;a,b)$ be the
$\{0,1,2,\ldots,\infty\}$-valued random variable giving the number of
\emph{distinct} points in $\R\times\{t_0+t\}$ that are touched by
paths in $V$ which also touch some point in $[a,b]\times \{t_0\}$.

The following is the characterization Theorem of the
Brownian web of Fontes, Isopi, Newman and Ravishankar in
\cite{1}.

{\Th A} {\it There is an $({\cal H},{\cal F}_{\cal H})$-valued
random variable $\bar \W$ whose distribution $\mu_{\bar W}$ is
uniquely determined by the following three properties:}

\vskip 2mm (o) {\it from any deterministic point $(x,t)$ in
$\R^2$, there is almost surely a unique path $W_{x,t}$ in 
$\bar \W$ starting from $(x,t)$.}

\vskip 2mm (i) {\it for any deterministic $n$ and
$(x_1,t_1),\ldots,(x_n,t_n)$, the joint
  distribution of $W_{x_1,t_1},\ldots$, $W_{x_n,t_n}$ is that of
  coalescing Brownian motions (with unit diffusion constant), and}

\vskip 2mm

(ii) {\it for any deterministic, dense countable subset ${\cal D}$
of $\R^2$,
  almost surely, $\bar \W$ is the closure in $({\cal H}, d_{\cal H})$ of
  $\{W_{x,t}:(x,t)\in {\cal D}\}$.}

\vskip 2mm
The $({\cal H},{\cal F}_{\cal H})$-valued random variable $\bar \W$ given in
Theorem A is called \emph{standard Brownian web}. Alternative
characterizations can be found in \cite{1} as well.

\vskip 2mm
\paragraph{Main result: invariance principle} 
The Poisson tree induces sets of continuous paths as follows. For any
$s=(s_1,s_2)\in S$, the Poisson process with parameter $\la$, we define the
path $X^s$ in $\R^2$ as the linearly interpolated line composed by all edges
$\{(\a^{n-1}(s),\a^n(s)):n\in\N\}$ of the Poisson tree $G$.  Clearly, $X^s\in
C[s_2]\times \{s_2\}\subset \Pi$. Let
\begin{equation}\label{100}X:=\{X^s:s\in S\},\end{equation} which we also call
the \emph{Poisson web}.

By the definition of the Poisson tree, $X$ depends on $\la>0$ and $r>0$. In
case of necessity, we denote it by $X(\la,r)$. Take $\la=\la_0={\sqrt 3}/6$,
$r=r_o=\sqrt 3$, and let
\begin{equation}
\label{10}X_1:=X(\la_0,r_0).
\end{equation}
\begin{equation}\label{11}
X_\delta:=\{(\delta x_1,\delta^2 x_2)\in \R^2:(x_1,x_2)\in
X_1\},\end{equation} 
for $\delta\in(0,1]$. Namely, $X_\delta$ is the
diffusive rescaling of $X_1$.

Another family of Poisson trees
$Y_\delta, \delta\in(0,1]$, is defined as
\begin{equation}\label{12}Y_\delta:=X(\la(\delta),r(\delta)),\end{equation}
with $\la(\delta)=\delta^{-1}$, $r(\delta)=({3\delta}/2)^{1/3}$ for
$\delta\in (0,1]$. It is straightforward to verify the following
lemma.

\begin{lemma}
  \label{1.1}
For any $\delta\in(0,1]$, almost surely, the
closures of the Poisson trees $X_\delta$ and $Y_\delta$ defined in
(\ref{10}), (\ref{11}) and (\ref{12}) are compact subsets of
($\Pi,d$).
\end{lemma}

By Lemma \ref{1.1}, the closure of $X_\delta$ (resp. $Y_\delta$), also denoted
by $X_\delta$ (resp. $Y_\delta$), which is
obtained by adding all the paths of the form $(f,t_0)$ with
$t_0\in [-\infty,\infty]$ and $f\equiv\infty$ or $f\equiv
-\infty$, is an $({\cal H},{\cal F}_{\cal H})$-valued random
variable.

Our main result is a proof that $X_\delta$ and $Y_\delta$
converge in distribution to the
Brownian web characterized in Theorem A. Comparing with the
classical Donsker's invariance principle \cite{7} for a single
path, we call it the \emph{invariance principle} in the web case.

{\Th 1.1} {\it Each of the rescaled Poisson trees $X_\delta$ and $Y_\delta$
  converges in distribution to the standard Brownian web as $\delta\rightarrow
  0$. } \vskip 2mm

\paragraph{Background} 
Arratia's construction of a system of coalescing one dimensional Brownian
motions starting from every space-time point in $\R^2$~\cite{4} begins with
coalescing Brownian motions starting from a countable dense subset of $\R^2$.
A kind of semicontinuity condition is then imposed in order to get paths
starting from every point of $\R^2$, and also to insure a certain flow
condition, which in particular yields a {\em unique} path starting at each
space-time point.  A variant of that construction was done by T\'oth and
Werner~\cite{5}, with a different kind of semicontinuity condition, also
yielding uniqueness, in work where the final object was auxiliary in the
definition and study of the {\em self-repelling motion}.  The construction of
Fontes, Isopi, Newman and Ravishankar~\cite{1,1a} was done with the aim of
providing a set up for weak convergence.  Their choice of $({\cal H},{\cal
  F}_{\cal H})$ as sample space, with its good topological properties, was
inspired by the set up of Aizenman and Burchard~\cite{6} for the study of
scaling limits of critical statistical mechanics models. This choice sets the
stage for the derivation of criteria for characterization of and weak
convergence to the constructed object, the Brownian web~\cite{1,1a}.  These
criteria were then verified for rescaled coalescing one dimensional random
walks starting from every space-time point in $\Z\times\R$~\cite{1,1a}.  The
main part of the proof of Theorem 1.1 consists of the verification of those
criteria for $X_\delta$ and $Y_\delta$. The construction of Fontes, Isopi,
Newman and Ravishankar begins as those of Arratia and of T\'oth and Werner,
with coalescing Brownian motions starting from a countable dense subset of
$\R^2$.  In order to have paths starting at every point, instead of imposing a
semicontinuity condition, in particular disregarding flow/uniqueness issues,
they take the closure in path space of the initial countable collection of
paths. The resulting set of paths, the Brownian web, turns out to be almost
surely compact, a property which allows it to live in $({\cal H},{\cal
  F}_{\cal H})$. Within the topological framework of that space, suitable
criteria for characterization of and weak convergence to the Brownian web are
then derived. The compactness of the Brownian web is its main distinguishing
feature vis-a-vis the previous constructions. One other such feature is the
occurrence in the Brownian web of {\em multiple} space-time points, where more
than one path start out from. Multiple points, with alternative definitions
but the same nature, are however a feature of all the constructions. The
semicontinuity conditions of Arratia and of T\'oth and Werner eliminate, each
in its own ad hoc way, all but one of the paths starting out from those
points. This elimination is at the root of the noncompactness of their
resulting sets of paths.

\vspace{.25cm}
We next state one of the above mentioned weak convergence criteria.
Let ${\cal D}$ be a countable dense set of points in $\R^2$.

{\Th B} \cite{1} {\it Suppose \( X_{1}, X_{2}, \dots \) are \(({\cal H}, {\cal
    F}_{{\cal H}}) \)-valued random variables with noncrossing paths. If, in
  addition, the following three conditions are valid, the distribution of
  $X_n$ converges to the distribution $\mu_{\bar{W}}$ of the standard Brownian
  web.}
\begin{itemize}
\item[$(I_1)$] {\it There exist $\theta^y_n\in X_n$ such that for any
    deterministic $y_1,\ldots,y_m\in{\cal D}$,
    ${\theta^{y_1}_n,\ldots,\theta^{y_m}_n}$ converge in distribution as $n\to
    \infty$ to coalescing Brownian motions (with unit diffusion constant)
    starting at} $y_1,\ldots,y_m$.
\item[$(B_1)$]\quad\,\,  $\limsup_{n \to \infty}\sup_{(a,t_0)\in\R^2}
\P(\eta_{_{X_n}}(t_0,t;a,a+\epsilon)\ge 2)\to0 \hbox{ as } \epsilon\to0+$;
\item[$(B_2)$]  $\epsilon^{-1}\limsup_{n \to \infty}\sup_{(a,t_0)\in\R^2}
\P(\eta_{_{X_n}}(t_0,t;a,a+\epsilon)\ge 3)\to0\hbox{ as } \epsilon\to0+$.
\end{itemize}

To prove the main result, we show in Section 2 
that the Poisson webs $X_\delta$ and
$Y_\delta$ satisfy the hypothesis of Theorem B.
The verification of $I_1$, on Subsection~\ref{findim}, relies on a 
comparison with independent paths and on the almost sure coalescence 
of the Poisson web paths with each other. See Lemma \ref{2.5}.

In Subsection~\ref{b1b2}, an FKG inequality enjoyed by the distribution 
of a single Poisson web path (Lemma~\ref{2.10}) and the $O(t^{-1/2})$ 
decay of the coalescence time of two such paths (Lemma~\ref{2.11}), 
combined with $I_1$, yield both $B_1$ and $B_2$. An analogous argument, which
relies as well on an FKG property of the constituent paths, also holds for
ordinary coalescing one dimensional random walks starting from all space time
points, and can be used, together with the analogue of Lemma~\ref{2.5}
(which follows immediately from Donsker's invariance principle in this case), 
to establish their convergence (when suitably rescaled) to the Brownian web.
The arguments for the Poisson web, which are similar in spirit
to the ones for coalescing random walks, are nonetheless much more involved
for the former case than the ones for the latter case --- these are essentially
immediate. 

Working out a second example of a process in the basin of attraction of the
Brownian web (the first example, just mentioned above, being ordinary one
dimensional coalescing random walks) that is natural on one side, and that
requires substantial technical attention on another side, is the primary point
of this paper. Its main result may have an applied interest, e.g.~in the
context of drainage networks. The convergence results here may lead to
rigorous/alternative verification of some of the scaling theory for those
networks. See~\cite{rr}.  Ordinary one dimensional coalescing random walks
starting from all space time points have also been proposed as model of a
drainage network~\cite{sch}, so the latter remark applies to them as well.
Another application would be in obtaining aging results from the scaling limit
results for systems that could be modelled by Poisson webs, like drainage
networks. For the relation between aging and scaling limits, see
e.g.~\cite{fins}, \cite{fin}, \cite{1}, \cite{1a} and references therein.


\section{Proofs}
\renewcommand{\theequation}{2.\arabic{equation}}
\setcounter{equation}{0}
\paragraph{Coalescing random walks} Let $S$ be the Poisson process
with parameter $\la>0$, fix some $r>0$. For any
$x=(x_1,x_2)\in\R^2$, let $\tau^n(x)=[\a^n(x)]_2$, $n\ge 0$, be the
second coordinate of $\a^n(x)$ and consider $\{\xi^x(t):t\ge x_2\}$ 
as the continuous time Markov process defined by
\begin{equation}\label{13}
 \xi^x(t)=[\a^n(x)]_1, {\rm the\ first\ coordinate\ of}\ \a^n(x);\
 t\in [\tau^n(x),\tau^{n+1}(x)),
 \ \ n\ge 0.
\end{equation} 
We remark that for any fixed $(x^i)_{i=1}^m$, with
$x^i=(x^i_1,x^i_2)\in\R^2$ for $i=1,\ldots,m$,
$\{(\xi^{x^i}(t):\,t\ge x^i_2),\, i=1,\ldots,m\}$ 
defines a finite system of coalescing random
walks starting at the space-time points $x^1,\ldots,x^m$.
For any fixed $x=(x_1,x_2)\in\R^2$, the marginal distribution 
of $\xi^x(\cdot)$ is that of a continuous 
time random walk starting at time $x_2$ in position $x_1$, which, at
exponentially distributed random waiting times of mean $(2r\la)^{-1}$, 
chooses a point uniformly in interval $[x_1-r,x_1+r]$ and then jumps 
to that point. 
The interaction appears when two walks are located at points 
$x_1\in\R$ and $y_1\in \R$ at some time $t_0$: if $|x_1-y_1|\le 2r$, 
and the Poisson process $S$ makes that
$\a(x_1,t_0)=\a(y_1,t_0)=s=(s_1,s_2)\in S$, then both walks jump
to the same position $s_1\in\R$ at time
$\tau(x_1,t_0)=\tau(y_1,t_0)=s_2$ and coalesce. We note that the
finite system of coalescing random walks 
$\{(\xi^{x^i}(t):\,t\ge x^i_2),\, i=1,\ldots,m\}$ is also strong
Markov.

For $x=(x_1,x_2)\in \R^2$, let
$x_\delta=(\delta^{-1}x_1,\delta^{-2}x_2)$, $\delta\in(0,1]$. For
the single random walk starting at $x=(x_1,x_2)$, $\xi^x(\cdot)$, defined 
in the last paragraph, the diffusive rescaling is
\begin{equation} \label{14}
\xi^x_\delta(t):=\delta\xi^{x_\delta}(\delta^{-2}t),\ \ \ {\rm
for}\ \; t\ge x_2; \ \ \ \delta\in (0,1].
\end{equation}

Since Theorems A and B apply to continuous paths only, we need to
replace the original processes by their linearly interpolated
versions:
\begin{equation}\label{15} \bar\xi^x_\delta(t)=
\delta\left\{\xi^{x_\delta}(\tau^n(x_\delta))+
\frac{\delta^{-2}t-\tau^n(x_\delta)}
{\tau^{n+1}(x_\delta)-\tau^n(x_\delta)}
\left(\xi^{x_\delta}(\tau^{n+1}(x_\delta))
-\xi^{x_\delta}(\tau^n(x_\delta))\right)\right\}, \end{equation} for $t\ge
x_2 $ such that
$\delta^{-2}t\in[\tau^n(x_\delta),\tau^{n+1}(x_\delta))$, $n\ge
0$; $\delta\in(0,1]$, $x\in \R^2$. Denote by $\bar\xi^x_\delta$
the corresponding continuous path in $\R^2$ and note that
$\bar\xi^s_1$ is just $X^s$ in (\ref{100}) with $s\in S$. 
It is straightforward to see that $\bar\xi^x_\delta\in X_\delta$, the
Poisson web defined by (\ref{11}), if and only if $x_\delta\in S$.

Let \begin{equation}\label{16}\theta^x_\delta:=\left\{\begin{array}{ll}
\bar\xi^x_\delta&,\ \ \ {\rm if}\ \ x_\delta\in S\\[3mm]
\bar\xi^{(\delta[\a(x_\delta)]_1,\delta^2[\a(x_\delta)]_2)}_\delta
&,\ \ \ {\rm otherwise.} \end{array}\right. \end{equation}
In this way, for all $x\in
\R^2$ and $\delta\in (0,1]$, $\theta^x_\delta\in X_\delta$. Note
that the paths defined by (\ref{15}) and (\ref{16}) depend on the
choice of $\la>0$ and $r>0$. In case of necessity,
we denote them by $\bar\xi^x_\delta(\la,r)$ and $\theta^x_\delta(\la,r)$.

\vskip 2mm
The following is an application of the classical Donsker's theorem
\cite{7} in our case.

\begin{lemma}
  \label{2.1}
If $\la=\la_0=\sqrt 3/6$, $r=r_0=\sqrt 3$, then
$\bar\xi^x_\delta$ converges in distribution as $\delta\rightarrow
0$ to $B^x$, the Brownian path with unit diffusion coefficient
starting from space-time point $x=(x_1,x_2)\in \R^2$.
\end{lemma}

\vskip 2mm For any $x^1,\ldots,x^m\in \R^2$, $m\in \N$, regard
$(\bar\xi^{x^1}_\delta,\ldots,\bar\xi^{x^m}_\delta)$ and
$(\theta^{x^1}_\delta,\ldots,\theta^{x^m}_\delta)$ as random variables in the
product metric space $(\Pi^m,d^{*m})$, where $d^{*m}$ is a distance on $\Pi^m$
such that the topology generated by it coincides with the corresponding
product topology. Here we choose and define \begin{equation}\label{17}
  d^{*m}[(\xi^1,\ldots,\xi^m),(\zeta^1,\ldots,\zeta^m)]=\max_{1\le i\le
    m}d(\xi^i,\zeta^i), \end{equation} for all
$(\xi^1,\ldots,\xi^m),(\zeta^1,\ldots,\zeta^m)\in\Pi^m$, where $d$ was defined
in (\ref{5}). The next result follows immediately from the definition.

\begin{lemma}
  \label{2.2}
\begin{equation}\label{18}
\P_\la\{d^{*m}[(\bar\xi^{x^1}_\delta,\ldots,\bar\xi^{x^m}_\delta),
(\theta^{x^1}_\delta,\ldots,\theta^{x^m}_\delta)]\ge
\epsilon\}\rightarrow 0,\ \ {\rm as}\ \ \delta\rightarrow 0\end{equation}
for all $\epsilon>0$, $\la>0$, $r>0$, and $x^1,\ldots,x^m\in \R^2,
m\in\N$, where $\P_\la$ is the probability distribution of $S$,
the Poisson process with parameter $\la$.
\end{lemma}

\vskip 2mm For the collection of Poisson trees $Y_\delta$'s, we
define the path $\bar\zeta^x_\delta$ as
\begin{equation}\label{19}\bar\zeta^x_\delta(t):=\bar\xi^x_1(t)(\la(\delta),r(\delta)),\
\ \ \forall \ \ x\in\R^2,\ \ \delta\in(0,1],  \end{equation} the rescaled
continuous path defined in (\ref{15}) with
$(\la,r)=(\la(\delta),r(\delta))=(\delta^{-1},(3\delta/2)^{1/3})$.
We also define the path \begin{equation}\label{20}
\vartheta^x_\delta:=\theta^x_1(\la(\delta),r(\delta)), \ \ \forall
\ \ x\in\R^2,\ \ \delta\in(0,1],  \end{equation} so that
$\vartheta^x_\delta\in Y_\delta$.

{}For the same reasons as in Lemma~\ref{2.1} and Lemma~\ref{2.2}, we have:
\begin{lemma}
  \label{2.3}
$\bar\zeta^x_\delta$ converges in distribution as
$\delta\rightarrow 0$ to $B^x$ for all space-time point
$x=(x_1,x_2)\in\R^2$.
\end{lemma}

\begin{lemma}
  \label{2.4}
\begin{equation}\label{21}
\P\{d^{*m}[(\bar\zeta^{x^1}_\delta,\ldots,\bar\zeta^{x^m}_\delta),
(\vartheta^{x^1}_\delta,\ldots,\vartheta^{x^m}_\delta)]\ge\epsilon\}
\rightarrow 0, \ {\rm as}\ \ \delta\rightarrow 0\end{equation} for all
$\epsilon>0$ and $x^1,\ldots,x^m\in \R^2, m\in\N$, where $d^{*m}$
is defined in (\ref{17}).
\end{lemma}


\vskip 2mm
\subsection{Convergence in finite-dimensional cases: verification of
condition $I_1$.}\label{findim} In this subsection, we begin to prove 
Theorem 1.1. In our proofs, we will mainly verify the corresponding
conditions of Theorem B for the Poisson trees $X_\delta,\delta\in(0,1]$,
because of the essential similarity with the Poisson trees 
$Y_\delta,\delta\in(0,1]$; some remarks will be given for the 
case of the latter processes.

Let ${\cal D}$ be a countable dense set of points in $\R^2$. 
As pointed in the last subsection,
for any $y\in \R^2$ and $\delta\in (0,1]$, as single-paths, 
$\theta^y_\delta\in X_\delta,\vartheta^y_\delta\in Y_\delta$. 
In this subsection, to prove condition $I_1$, we will show that, for any
$y^1,y^2,\ldots,y^m\in {\cal D}$,
$(\theta^{y^1}_\delta,\ldots,\theta^{y^m}_\delta)$ and
$(\vartheta^{y^1}_\delta,\ldots,\vartheta^{y^m}_\delta)$ converge
in distribution as $\delta\rightarrow 0$ to coalescing Brownian
motions (with unit diffusion constant) starting at
$y^1,\ldots,y^m$. Actually, by Lemma~\ref{2.2} and Lemma~\ref{2.4}, we only
need to prove the following.

\begin{lemma}
  \label{2.5}
$(\bar\xi^{y^1}_\delta,\ldots,\bar\xi^{y^m}_\delta)$ and
$(\bar\zeta^{y^1}_\delta,\ldots,\bar\zeta^{y^m}_\delta)$ converge
in distribution as $\delta\rightarrow 0$ to coalescing Brownian
motions (with unit diffusion constant) starting at
$y^1,\ldots,y^m$ ($\in {\cal D}$).
\end{lemma}

For the finite system of coalescing random walks defined in the
last subsection, Ferrari, Landim and Thorisson \cite{2}
proved that, for any $x^1,x^2\in\R^2$, the random walks $\xi^{x^1}(t)$
and $\xi^{x^2}(t)$, $t\ge x^1_2\vee x^2_2$ will meet and then
coalesce almost surely. This also follows from Lemma 2.10 below.
The following is a corollary of this result.

\begin{lemma}
  \label{2.6}
For any $\la>0$, $r>0$, we have \begin{equation}\label{22}
\lim_{\sigma\rightarrow\infty}\P_\la\{\bar
d(\bar\xi^{x^1}_1,\bar\xi^{x^2}_1)\ge\sigma\}=0, \end{equation} for all
$x^1,x^2\in \R^2$, where $\bar d$ is a function defined on $\Pi^2$
such that for any $(f_1,t_1),(f_2,t_2)\in\Pi$ \begin{equation}\label{22'} \bar
d((f_1,t_1),(f_2,t_2)):=\sup_t|\hat f_1(t)-\hat
f_2(t)|\vee|t_1-t_2|. \end{equation}
\end{lemma}
 
Now, for any $y^1,\ldots,y^m\in{\cal D}$, let
$B^{y^1},\ldots,B^{y^m}$ be $m$ independent Brownian paths
starting at space-time points $y^1,\ldots,y^m$, respectively. As
Arratia did in \cite{3}, we construct the one-dimensional
coalescing Brownian motions starting at $y^1,\ldots,y^m$ by
defining a continuous function $f$ from $\Pi^m$ to $\Pi^m$.

Let $\gamma_0=\min(y^1_2,\ldots,y^m_2)$, define the stopping time
\begin{equation}\label{t1} \gamma_1=\min\{t>\gamma_0:\exists \ 1\le i,j\le m,
{\rm such\ that}\ \ t\ge y^i_2\vee y^j_2, \ \ {\rm and}\ \
B^{y^i}(t)=B^{y^j}(t)\}. \end{equation}

At time $\gamma_1$, let $(i,j)$ be the pair such that
 $B^{y^i}(\gamma_1)=B^{y^j}(\gamma_1)$ and $B^{y^i}(t)<B^{y^j}(t)$
for all $t\in [y^i_2\vee y^j_2,\gamma_1)$, note that such pair
$(i,j)$ is unique almost surely and we denote ${\cal I}_1=\{i\}$,
${\cal J}_1=\{j\}$. Denote by $p_1:=B^{y^{i}}(\gamma_1)$ the
position of the first coalescence. Now, we renew the system by
resetting $B^{y^j}(t)=B^{y^i}(t)$ for all $t\ge \gamma_1$.

For all $1\le k\le m-2$, after we have defined
$\gamma_1,\ldots,\gamma_k$ and renewed the system $k$ times,
we define \begin{equation}\label{tk}
\gamma_{k+1}:=\min\left\{t>\gamma_{k}:\hskip -3mm \begin{array}{ll}&
\exists\, 1\le i,j\le m,\ {\rm such\ that}\ \ t\ge y^i_2\vee
y^j_2, \ {\rm and}\
B^{y^i}(t)\\[2mm]&=B^{y^j}(t),\ B^{y^i}(t')\not= B^{y^j}(t'),\ \forall
\ t'\in[y^i_2\vee y^j_2,t).\end{array} \right\}. \end{equation} Let $(i,j)$
be the (almost surely) unique pair such that
$B^{y^i}(\gamma_{k+1})=B^{y^j}(\gamma_{k+1})$, and
$B^{y^i}(t)<B^{y^j}(t)$ for all $t\in[y^i_2\vee
y^j_2,\gamma_{k+1})$. Let \[{\cal I}_{k+1}=\{i':\exists\
\epsilon>0,\ {\rm such}\ {\rm that}\ B^{y^{i'}}(t)=B^{y^i}(t),\
\forall\ t\in[\gamma_{k+1}-\epsilon,\gamma_{k+1})\};\]
\[{\cal
J}_{k+1}=\{j':\exists\ \epsilon>0,\ {\rm such}\ {\rm that}\
B^{y^{j'}}(t)=B^{y^j}(t),\ \forall\
t\in[\gamma_{k+1}-\epsilon,\gamma_{k+1})\}.
\] We renew the system for the $k+1$-st time by resetting
$B^{y^j}(t)=B^{x^{i_{k+1}}}(t)$ for all $j\in {\cal I
}_{k+1}\cup{\cal J}_{k+1}\setminus\{i_{k+1}\}$ and
$t\ge\gamma_{k+1}$, where $i_{k+1}\in{\cal I}_{k+1}$ satisfying
that, for any other $i'\in{\cal I}_{k+1}$, $B^{y^{i'}}(t)\ge
B^{y^{i_{k+1}}}(t)$ for all $t\in [y^{i_{k+1}}_2\vee
y^{i'}_2,\gamma_{k+1})$ . Denote by
$p_{k+1}:=B^{y^{i_{k+1}}}(\gamma_{k+1})$ the position of the
$k+1$-st coalescence.

By the basic properties of one-dimensional Brownian motion, we
have that \begin{equation}\label{23} \begin{array}{ll}
&-\infty<\gamma_0<\gamma_1<\ldots<\gamma_{m-1}<\infty\ {\rm
and}\\[2mm]
& \varrho:=\min\{|p_k-p_{k'}|: 1\le k,k'\le m-1\}>0 \end{array} \end{equation} almost
surely. That is, let $C\subset\Pi^m$ be the set of all $m$-dimensional
continuous paths starting at space-time points $y^1,\ldots,y^m\in
\R^2$ satisfying condition (\ref{23}). Then, we have \begin{equation}\label{dc}
\P\{(B^{y^1},\ldots,B^{y^m})\in C\}=1. \end{equation}

The resulting system after $m-1$ steps of renewing is the so-called
one-dimensional coalescing Brownian motions starting at space-time
points $y^1,\ldots,y^m$, which is denoted by
$f(B^{y^{1}},\ldots,B^{y^{m}})$, a function of the $m$ independent
Brownian motions $B^{y^1},\ldots,B^{y^m}$.

For any $m$ distinct points 
$y^1,\ldots,y^m\in{\cal D}$ and $\delta
\in(0,1]$. Let $\bar\xi^{y^1}_\delta,\ldots,\bar\xi^{y^m}_\delta$
be the $m$ rescaled continuous random paths defined in (\ref{15})
from the same Poisson process with $\la=\la_0=\sqrt 3/6$ and
$r=r_0=\sqrt 3$. Having
$(\bar\xi^{y^1}_\delta,\ldots,\bar\xi^{y^m}_\delta)$ as a random
element in $\Pi^m$, we want to define a function $f_\delta$ of it
to $\Pi^m$. This is our main idea for the
verification of condition $I_1$: we define what we call
``$\delta$-coalescence" of the random paths
$\bar\xi^{y^1}_\delta,\ldots,\bar\xi^{y^m}_\delta$ 
in such a way that, in the system
$f_\delta(\bar\xi^{y^1}_\delta,\ldots,\bar\xi^{y^m}_\delta)$,
before any $\delta$-coalescence, the paths involved are
independent.

Similarly as we have done with $f$ in the preceding paragraphs, we 
define $f_\delta$ by renewing the whole system step by step as 
follows.
Let $\gamma_{\delta,0}=\min(y^1_2,\ldots,y^m_2)$. From now on, we 
assume that $\delta>0$ is close enough to $0$, so that, in particular,
the following stopping time is well-defined. 
\be
\label{24}
\gamma_{\delta,1}=\inf\{t>\gamma_{\delta,0}:\exists 1\le i,j\le m,
{\rm such\ that}\ t\ge y^i_2\vee y^j_2,\
|\xi^{y^i}_\delta(t)-\xi^{y^j}_\delta(t)|<2\sqrt 3\delta\}, 
\end{equation}
where $\xi^{y^i}_\delta(t), 1\le i\le m$ is the rescaled random
walk defined in (\ref{14}).

Suppose that $(i,j)$ is the (almost surely) unique pair such that
$|\xi^{y^i}_{\delta}(\gamma_{\delta,1})-\xi^{y^j}_\delta(\gamma_{\delta,1})|\le
2 \sqrt 3\delta$ and $\xi^{y^j}_{\delta}(t)-\xi^{y^i}_\delta(t)>
2\sqrt 3\delta$ for all $t\in[y^i_2\vee y^j_2,\gamma_{\delta,1})$.
We denote it by $(i_1,j_1)$, and let ${\cal
I}_{\delta,1}:=\{i_1\}$, ${\cal J}_{\delta,1}:=\{j_1\}$. We renew
the system according to the following two cases. Case (a): There
exists $n'\ge 1$ such that
$\delta^{-2}\gamma_{\delta,1}=\tau^{n'}(y^{j_1}_\delta)$. Case
(b): There exists $n^{\prime\prime}\ge 1$ such that
$\delta^{-2}\gamma_{\delta,1}=\tau^{n''}(y^{i_1}_\delta)$. Recall
that $y^{i_1}_\delta=(\delta^{-1}
y^{i_1}_1,\delta^{-2}y^{i_1}_2)\in \R^2$, i.e., the starting
space-time point of the original (before rescaling) random walk of
$\xi^{y^{i_1}}_\delta$.

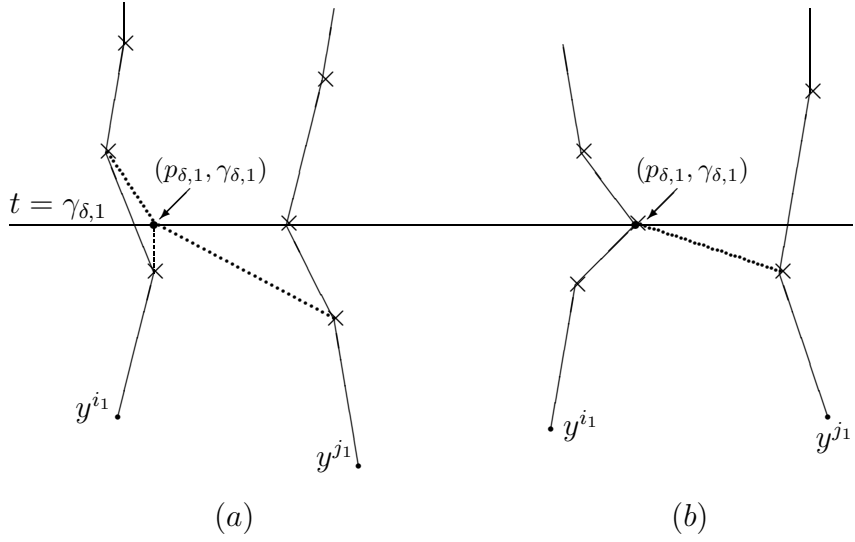
\begin{figure}
\unitlength=0.8 mm
\begin{picture}(120,70)(-24,-5)
\put(18,8){\circle*{1}}\put(18,8){\line(1,4){6}}\put(22.2,31){$\times$}\put(11,8){$y^{i_1}$}
\put(24,32){\line(-2,5){8}}\put(14.4,51){$\times$}
\put(16,52){\line(1,6){3}}\put(17.3,69){$\times$}
\put(51,0){$y^{j_1}$}\put(58,0){\circle*{1}}\put(58,0){\line(-1,6){4}}\put(19,70){\line(0,1){7}}
\put(52.2,23.2){$\times$}\put(54,24){\line(-1,2){8}}
\put(46,40){\line(1,4){6}}\put(44.4,39){$\times$}\put(52,64){\line(1,6){2}}
\put(50.5,63){$\times$} \put(0,40){\line(1,0){140}}
\put(0,42){$t=\gamma_{\delta,1}$}
\multiput(24,32)(0,1){8}{\line(0,1){0.5}}
\multiput(23.6,38.8)(-0.5,0.75){17}{$\cdot$}
\multiput(53.2,23.2)(-0.9385,0.5){32}{$\cdot$}
\put(24,40){\circle*{1.5}} \put(34,-10){($a$)}\put(110,-10){($b$)}
\put(24,48){\footnotesize($p_{\delta,1},\gamma_{\delta,1}$)}
\put(30,46){\vector(-1,-1){5}}
\put(90,6){\line(1,6){4}}\put(92.4,28.9){$\times$}
\put(94,30){\line(1,1){10}} \put(104,40){\line(-3,4){9}}
\put(95,52){\line(-1,6){3}}\put(93.5,51){$\times$}
\put(136,8){\circle*{1}}\put(136,8){\line(-1,3){8}}\put(126.5,31){$\times$}
\multiput(127.4,31)(-0.75,0.25){33}{$\cdot$}
\put(128,32){\line(1,6){5}}\put(133,62){\line(0,1){14}}\put(131.5,61){$\times$}
\put(104,40){\circle*{1.5}}
\put(104,48){\footnotesize($p_{\delta,1},\gamma_{\delta,1}$)}
\put(111,46){\vector(-1,-1){5}}
\put(90,6){\circle*{1}}\put(92,5){$y^{i_1}$}
\put(134,3){$y^{j_1}$}\put(102.4,39){$\times$}
\end{picture}
\vskip 10mm
\begin{center}
\begin{minipage}{12 cm}
\caption{\normalsize  The points marked by ``$\times$" are those
space-time points $x$'s such that $x_\delta\in S$, the points of
the Poisson process. To renew the system, one may go along with
the dash lines if possible. In both cases (a) and (b), in the
renewed system, $\bar\xi^{y^{i_1}}_\delta$ and
$\bar\xi^{y^{j_1}}_\delta$ meet and coalesce at space-time point
($p_{\delta,1},\gamma_{\delta,1}$). }
\end{minipage}
\end{center}
\end{figure}

In case (a), suppose that
$\delta^{-2}\gamma_{\delta,1}=\tau^{n'}(y^{j_1}_\delta)\in
[\tau^{k'}(y^{i_1}_\delta),\tau^{k'+1}(y^{i_1}_\delta))$ for some
$k'\ge 0$. We renew the system by resetting
$\bar\xi^{y^{i_1}}_\delta(t)$ for $t\in
[\delta^2\tau^{k'}(y^{i_1}_\delta),\delta^2\tau^{k'+1}(y^{i_1}_\delta))$
and resetting $\bar \xi^{y^{j_1}}_\delta(t)$ for all
$t\ge\delta^2\tau^{n'-1}(y^{j_1}_\delta)$ as follows (see Figure
1(a)):

\begin{equation}\label{25} \bar\xi^{y^{i_1}}_\delta(t)=\left\{\begin{array}{ll}
\delta\cdot\xi^{y^{i_1}_\delta}(\tau^{k'}(y^{i_1}_\delta))&\ \ \
,\ \ \
t\in[\delta^{2}\tau^{k'}(y^{i_1}_\delta),\gamma_{\delta,1})\\[8mm]
\delta\cdot\xi^{y^{i_1}_\delta}(\tau^{k'}(y^{i_1}_\delta))+
\delta\cdot\frac{\delta^{-2}t-\tau^{n'}(y^{j_1}_\delta)}
{\tau^{k'+1}(y^{i_1}_\delta)-\tau^{n'}(y^{j_1}_\delta)}&\\[-1mm]
&\ \ \  ,\ \ \
t\in[\gamma_{\delta,1},\delta^2\tau^{k'+1}(y^{i_1}_\delta))\\[-1mm]
\cdot\left[\xi^{y^{i_1}_\delta}(\tau^{k'+1}(y^{i_1}_\delta))
-\xi^{y^{i_1}_\delta}(\tau^{k'}(y^{i_1}_\delta))\right] & \end{array}
\right.\hskip3mm \end{equation}

\vskip 3mm \begin{equation}\label{26}
\bar\xi^{y^{j_1}}_\delta(t)=\left\{\begin{array}{ll}
\bar\xi^{y^{i_1}}_\delta(t)&,\ \ \ t\ge\gamma_{\delta,1}\\[8mm]
\delta
\cdot\xi^{y^{j_1}_\delta}(\tau^{n'-1}(y^{j_1}_\delta))+\delta\cdot\frac{\delta^{-2}t-\tau^{n'-1}(y^{j_1}_\delta)}
{\tau^{n'}(y^{j_1}_\delta)-\tau^{n'-1}(y^{j_1}_\delta)}&\\[-1mm]
&,\ \ \
t\in[\delta^{2}\tau^{n'-1}(y^{j_1}_\delta),\gamma_{\delta,1})\\[-1mm]
\cdot\left[\xi^{y^{i_1}_\delta}(\tau^{k'}(y^{i_1}_\delta))
-\xi^{y^{j_1}_\delta}(\tau^{n'-1}(y^{j_1}_\delta))\right]& \end{array}
\right.\hskip 1.5mm \end{equation} In the renewed system, the paths
$\bar\xi^{y^{i_1}}_\delta$ and $\bar\xi^{y^{j_1}}_\delta$ meet and
coalesce at time $\gamma_{\delta,1}$.

In case (b), suppose that
$\delta^{-2}\gamma_{\delta,1}=\tau^{n''}(y^{i_1}_\delta)\in
[\tau^{k''}(y^{j_1}_\delta),\tau^{k''+1}(y^{j_1}_\delta))$ for
some $k''\ge 0$. We renew the system by resetting $\bar
\xi^{y^{j_1}}_\delta(t)$ for all
$t\ge\delta^2\tau^{k''}(y^{j_1}_\delta)$ as follow (see Figure 1(b)).

\begin{equation}\label{27} \bar\xi^{y^{j_1}}_\delta(t)=\left\{\begin{array}{ll}
\bar\xi^{y^{i_1}}_\delta(t)&\ \ \ \ ,\ \ \ t\ge\gamma_{\delta,1}\\[8mm]
\delta
\cdot\xi^{y^{j_1}_\delta}(\tau^{k''}(y^{j_1}_\delta))+\delta\cdot\frac{\delta^{-2}t-\tau^{k''}(y^{j_1}_\delta)}
{\tau^{n''}(y^{i_1}_\delta)-\tau^{k''}(y^{j_1}_\delta)}&\\[-1mm]
&\ \ \ \ ,\ \ \
t\in[\delta^{2}\tau^{k''}(y^{j_1}_\delta),\gamma_{\delta,1})\\[-1mm]
\cdot\left[\xi^{y^{i_1}_\delta}(\tau^{n''}(y^{i_1}_\delta))
-\xi^{y^{j_1}_\delta}(\tau^{k''}(y^{j_1}_\delta))\right]& \end{array}
\right. \hskip 5mm\end{equation} Denote by
$p_{\delta,1}:=\xi^{y^{i_1}}_\delta(\gamma_{\delta,1})$ the
position of the first $\delta$-coalescence, and finally let ${\cal
K}_1=\{1,2,\ldots,m\}\setminus\{j_1\}$.

For all $1\le k\le m-2$, after we have defined
$\gamma_{\delta,1},\ldots,\gamma_{\delta,k}$; ${\cal
K}_1,\ldots,{\cal K}_k$ and renewed the system $k$ times,
we define \begin{equation}\label{28}
\gamma_{\delta,k+1}:=\inf\left\{t>\gamma_{\delta,k}:\hskip -3mm
\begin{array}{ll}& \exists\ i,j\in {\cal K}_k,\ {\rm such\ that}\ t\ge
y^i_2\vee
y^j_2, \ |\xi^{y^i}_\delta(t)-\xi^{y^j}_\delta(t)|<\\[2mm]& 2\sqrt 3\delta\
{\rm and}\ \xi^{y^j}_\delta(t')- \xi^{y^i}_\delta(t')>2\sqrt
3\delta,\ \forall \ t'\in[y^i_2\vee y^j_2,t)\end{array} \right\}. \end{equation} Let
$(i,j)$ be the (almost surely) unique pair such that
$|\xi^{y^i}_\delta(\gamma_{\delta,k+1})-\xi^{y^j}_\delta
(\gamma_{\delta,k+1})|<2\sqrt 3\delta$ and 
$\xi^{y^j}_\delta(t)- \xi^{y^i}_\delta(t)>2\sqrt 3\delta$ 
for all $t\in[\gamma_{\delta,k},\gamma_{\delta,k+1})$.
Let \[{\cal I}_{\delta,k+1}=\{i':\exists\ \epsilon>0,\ {\rm such}\
{\rm that}\ \bar\xi^{y^{i'}}_\delta (t)=\bar\xi^{y^i}(t),\
\forall\
t\in[\gamma_{\delta,k+1}-\epsilon,\gamma_{\delta,k+1})\};\]
\[{\cal
J}_{\delta,k+1}=\{j':\exists\ \epsilon>0,\ {\rm such}\ {\rm that}\
\bar\xi^{y^{j'}}_\delta(t)=\bar\xi^{y^j}_\delta(t),\ \forall\
t\in[\gamma_{\delta,k+1}-\epsilon,\gamma_{\delta,k+1})\}.
\]
Let $i_{k+1}\in{\cal I}_{\delta,k+1}$ (resp. $j_{k+1}\in{\cal
J}_{\delta,k+1}$) satisfy that, for any other $i'\in{\cal
I}_{\delta,k+1}$ (resp. $j'\in{\cal J}_{\delta,k+1}$),
$\bar\xi^{y^{i'}}(t)\ge \bar\xi^{y^{i_{k+1}}}(t)$ (resp.
$\bar\xi^{y^{j'}}(t)\ge \bar\xi^{y^{j_{k+1}}}(t)$) for all $t\in
[y^{i_{k+1}}_2\vee y^{i'}_2,\gamma_{\delta,k+1})$ (resp. $t\in
[y^{j_{k+1}}_2\vee y^{j'}_2,\gamma_{\delta,k+1})$).

Now, suppose that at time $\gamma_{\delta,k_0}$, $1\le k_0\le k$
(resp. $\gamma_{\delta,k_1}$, $1\le k_1\le k$), all paths of
$\{\bar\xi^{y^i}_\delta:i\in{\cal I}_{\delta,k+1}\}$ (resp.
$\{\bar\xi^{y^j}_\delta:j\in{\cal J}_{\delta,k+1}\}$) have coalesced
into one path. Then we renew the system in two steps. Firstly, we
renew $\bar\xi^{y^{i_{k+1}}}_\delta$ and
$\bar\xi^{y^{j_{k+1}}}_\delta$. Because it is only necessary to
renew $\bar\xi^{y^{i_{k+1}}}_\delta$ and
$\bar\xi^{y^{j_{k+1}}}_\delta$ after time $\gamma_{\delta,k_0}$
and $\gamma_{\delta,k_1}$, respectively, we deal with that as we did
in equations (\ref{25}), (\ref{26}) and (\ref{27}) for the
two-path system
$(\bar\xi^{(p_{\delta,k_0},\gamma_{\delta,k_0})}_\delta,
\bar\xi^{(p_{\delta,k_1},\gamma_{\delta,k_1})}_\delta)$
starting from space-time points
$(p_{\delta,k_0},\gamma_{\delta,k_0}),
(p_{\delta,k_1},\gamma_{\delta,k_1})\in\R^2$.
Secondly, we renew the system by resetting
$\bar\xi^{y^i}_\delta(t)=\bar\xi^{y^{i_{k+1}}}_\delta(t)$ for all
$i\in{\cal I}_{\delta,k+1}$ and $t\ge\gamma_{\delta,k_0}$, and
resetting
$\bar\xi^{y^j}_\delta(t)=\bar\xi^{y^{j_{k+1}}}_\delta(t)$ for all
$j\in{\cal J}_{\delta,k+1}$ and $t\ge\gamma_{\delta,k_1}$. Let
${\cal K}_{k+1}={\cal K}_{k}\setminus\{j_{k+1}\}$, and
$p_{\delta,k+1}:=\xi^{y^{i_{k+1}}}_\delta(\gamma_{\delta,k+1})$,
the position of the $(k+1)$th $\delta$-coalescence.

We denote the resulting object after renewing $m-1$ times by
$f_\delta(\bar\xi^{y^1}_\delta,\ldots,\bar\xi^{y^m}_\delta)$ and,
with that,
finish the definition of the function $f_\delta$. Clearly, for 
$\delta\in (0,1]$ small enough, we have
\be
\label{gadek}
-\infty<\gamma_{\delta,0}<\gamma_{\delta,1}<\ldots<\gamma_{\delta,m-1}<\infty
\end{equation} 
and the function $f_\delta$ is well defined almost surely.

Now, suppose that $\tilde\xi^{y^i}_\delta$ has the same
distribution as $\bar\xi^{y^i}$, $1\le i\le m$, and, as a random
element in $\Pi^m$,
$(\tilde\xi^{y^1}_\delta,\ldots,\tilde\xi^{y^m}_\delta)$ has
independent components. It is easy to see that the function $f_\delta$
is also well defined for the random paths
$(\tilde\xi^{y^1}_\delta,\ldots,\tilde\xi^{y^m}_\delta)$. Let
$C_\delta\subset\Pi^m$ be such that
\begin{equation}\label{dcd}\P\{(\tilde\xi^{y^1}_\delta,\ldots,\tilde\xi^{y^m}_\delta)\in
C_\delta\}=1\end{equation} and, on $C_\delta$, $f_\delta$ is well defined.

\begin{lemma}
  \label{2.7} Let
$(\bar\xi^{y^1}_\delta,\ldots,\bar\xi^{y^m}_\delta)$ is the $m$
rescaled continuous random paths defined in (\ref{15}) from the
same Poisson process with $\la=\la_0=\sqrt 3/6$ and $r=r_0=\sqrt
3$, and $(\tilde\xi^{y^1}_\delta,\ldots,\tilde\xi^{y^m}_\delta)$
have independent components and $\tilde\xi^{y^i}_\delta$ have the
same distribution as $\bar\xi^{y^i}$ for all $1\le i\le m$. Then,
\vskip -2mm 
\begin{description}
\vskip -2mm 
\item[\ \ \ \ \ a)]
$f_\delta(\bar\xi^{y^1}_\delta,\ldots,\bar\xi^{y^m}_\delta)$ has
the same distribution as
$f_\delta(\tilde\xi^{y^1}_\delta,\ldots,\tilde\xi^{y^m}_\delta)$.
\item[\ \ \ \ \ b)]
$f_\delta(\tilde\xi^{y^1}_\delta,\ldots,\tilde\xi^{y^m}_\delta)$
converges in distribution to $f(B^{y^1},\ldots,B^{y^m})$ as
$\delta\rightarrow 0$.
\item[\ \ \ \ \ c)] for any $\epsilon>0$,
\vskip -5.5mm 
\begin{equation}\label{29} \hskip 5mm
\P\{d^{*m}[f_\delta(\bar\xi^{y^1}_\delta,\ldots,\bar\xi^{y^m}_\delta),
(\bar\xi^{y^1}_\delta,\ldots,\bar\xi^{y^m}_\delta)]\ge\epsilon\}\rightarrow
0,
\ \ {\rm as}\ \delta\rightarrow 0, \end{equation} where $d^{*m}$ was
defined in (\ref{17}).
\end{description}
\end{lemma}

\vskip 2mm \proof a) Immediate from the definition of $f_\delta$. 
We only need to prove b) and c).

By Lemma~\ref{2.1} and independence,
$(\tilde\xi^{y^1}_\delta,\ldots,\tilde\xi^{y^m}_\delta)$ converges
in distribution to $(B^{y^1},\ldots$, $B^{y^m})$ as
$\delta\rightarrow 0$. By an extended continuous mapping theorem
of Mann and Wald \cite{8}, Prohorov \cite{9} (see also Theorem
3.27 of \cite{10}), we only need to prove that, for any
$c=(c^1,\ldots,c^m)\in C$, if
$c_\delta=(c^1_\delta,\ldots,c^m_\delta)\in C_\delta $ such that
$d^{*m}(c_\delta,c)\rightarrow 0$ as $\delta\rightarrow 0$, then
$d^{*m}(f_\delta(c_\delta),f(c))\rightarrow 0$ as
$\delta\rightarrow 0$. It is straightforward to check this by the
definitions of $f_\delta$, $f$, and $d^{*m}$; then we get b).

Since the function $\tanh(t)$ is Lipschitz continuous, we have
\begin{equation}\label{29'}
d^{*m}[f_\delta(\bar\xi^{y^1}_\delta,\ldots,\bar\xi^{y^m}_\delta),
(\bar\xi^{y^1}_\delta,\ldots,\bar\xi^{y^m}_\delta)]\le C_L\bar
d^{*m}
[f_\delta(\bar\xi^{y^1}_\delta,\ldots,\bar\xi^{y^m}_\delta),
(\bar\xi^{y^1}_\delta,\ldots,\bar\xi^{y^m}_\delta)],\end{equation} where
$0<C_L<\infty$ is the corresponding Lipschitz constant, and $\bar
d^{*m}$ is defined from $\bar d$ in the same way as we did in
(\ref{17}).

Now, if $\bar d^{*m}
[f_\delta(\bar\xi^{y^1}_\delta,\ldots,\bar\xi^{y^m}_\delta),
(\bar\xi^{y^1}_\delta,\ldots,\bar\xi^{y^m}_\delta)]\ge C_L\cdot
\epsilon$, then by the definition of function $f_\delta$, there
should exist some $1\le k\le m-1$ such that \begin{equation}\label{29"} \bar
d(\bar\xi^{(p_{\delta,k},\gamma_{\delta,k})}_\delta,
\bar\xi^{(\xi^{y^{j_k}}_\delta(\gamma_{\delta,k}),
\gamma_{\delta,k})}_\delta)\ge C_L\cdot \epsilon/n(j_k), \end{equation} where
$n(j_k)$ is the number of times that path
$\bar\xi^{y^{j_k}}_\delta$ was finally renewed in the process of
the definition of $f_\delta$. Obviously, $1\le n(j_k)\le m-1$.

By the basic property of coalescence, we have
\begin{equation}\label{30'} \bar
d(\bar\xi^{(p_{\delta,k},\gamma_{\delta,k})}_\delta,
\bar\xi^{(\xi^{y^{j_k}}_\delta(\gamma_{\delta,k}),
\gamma_{\delta,k})}_\delta)\le \bar
d(\bar\xi^{(p_{\delta,k},\gamma_{\delta,k})}_\delta,
\bar\xi^{(p_{\delta,k}+2\sqrt 3\delta,
\gamma_{\delta,k})}_\delta).\end{equation}

Then, by (\ref{29'}), (\ref{29"}), (\ref{30'}), the strong
Markov property of the coalescing random walks, and also the stationarity
of the Poisson process, we have \begin{equation}\label{30"} \begin{array}{ll}
&\P\{d^{*m}[f_\delta(\bar\xi^{y^1}_\delta,\ldots,\bar\xi^{y^m}_\delta),
(\bar\xi^{y^1}_\delta,\ldots,\bar\xi^{y^m}_\delta)]\ge\epsilon\}\\[3mm]
&\le\P\{\bar
d^{*m}[f_\delta(\bar\xi^{y^1}_\delta,\ldots,\bar\xi^{y^m}_\delta),
(\bar\xi^{y^1}_\delta,\ldots,\bar\xi^{y^m}_\delta)]\ge
C_L\cdot\epsilon\}\\[3mm]
&\le \Dp\sum_{k=1}^{m-1}\P\{\bar
d(\bar\xi^{(p_{\delta,k},\gamma_{\delta,k})}_\delta,
\bar\xi^{(p_{\delta,k}+2\sqrt 3\delta,
\gamma_{\delta,k})}_\delta)\ge C_L\cdot\epsilon/n(j_k)\}\\[3mm]
&\le(m-1)\P\{\bar d(\bar\xi^{(0,0)}_\delta,
\bar\xi^{(2\sqrt 3\delta,0)}_\delta)\ge
C_L\cdot\epsilon/(m-1)\}\\[3mm]
&=(m-1)\P_{\la_0}\{\bar d(\bar\xi^{(0,0)}_1,\bar\xi^{(2\sqrt
3,0)}_1)\ge {\delta^{-1}\cdot C_L\cdot \epsilon}/{(m-1)}\}. \end{array}
\end{equation} Now, by Lemma~\ref{2.6}, we get c), and the proof is finished.
$\square$

Lemma~\ref{2.5} is an immediate consequence of Lemma~\ref{2.7}. Thus,
condition $I_1$ for the Poisson 
web $X_\delta, \delta\in(0,1]$ follows from Lemma~\ref{2.2}.


{\Re To verify condition $I_1$ for the Poisson trees
$Y_\delta,\delta\in(0,1]$, in the definition of the function
$f_\delta$, one should use ``$\delta$-coalescence" when the
distance of two rescaled random walks is less than
$2r(\delta)=2(3\delta/2)^{1/3}$. Recall that for the Poisson trees
$X_\delta,\delta\in(0,1]$, we use it when that distance
is less than $2\sqrt 3\delta$.}


\subsection{Verification of conditions $B_1$ and $B_2$}\label{b1b2}
Consider the Poisson process $S$ with parameter $\la>0$ and
the corresponding Poisson tree $X:=X(\la, r)$ defined in
(\ref{100}) with respect to some fixed $r>0$. Given $t_0\in \R$, $t>0$,
$a,b\in\R$ with $a<b$, let $\eta_{_X}(t_0,t;a,b)$ be the
$\{0,1,2,\ldots,\infty\}$-valued random variable defined before
the statement of Theorem A. Let $\bar\eta_{_X}(t_0,t;a,b)$ be
another $\{0,1,2,\ldots,\infty\}$-valued random variable defined as the
number of distinct points $y=(y_1,y_2)\in
\R\times\{t_0+t\}$ such that there exists $s\in S$ with $s_2\le
t_0$, $\xi^s(t_0)\in [a,b]$ and $\xi^s(t_0+t)=\xi^s(y_2)=y_1$,
where $\xi^s$ is the Markov process defined in (\ref{13}). 
It is straightforward to see that, for any fixed $n\in\N$,
\begin{equation}\label{1000} \bar\eta_{_X}(t_0,t;a,b)\ge n\Rightarrow
\eta_{_X}(t_0,t;a-2r,b+2r)\ge n\Rightarrow
\bar\eta_{_X}(t_0,t;a-4r,b+4r)\ge n.\end{equation} This implies that, to
verify conditions $B_1$, $B_2$ for Poisson trees $X_\delta$ and
$Y_\delta$, we only need to verify the following $B_1'$ and
$B_2'$ respectively.

\begin{itemize}
\item[$(B_1')$]\quad\,\,$\limsup_{n \to \infty}
\P(\bar\eta_{_{\delta_n}}(0,t;0,\epsilon)\ge 2)\to0 \hbox{ as } \epsilon\to0+$;
\item[$(B_2')$]\quad  $\epsilon^{-1}\limsup_{n \to \infty}
\P(\bar\eta_{_{\delta_n}}(0,t;0,\epsilon)\ge 3)\to0\hbox{ as } \epsilon\to0+$
\end{itemize}
for any sequence of positive numbers $(\delta_n)$ such that 
$\lim_{n\to\infty}\delta_n=0$,
where
$\bar\eta_\delta=\bar\eta_{_{X_\delta}}$ or
$\bar\eta_{_{Y_\delta}}$, and we have used the space homogeneity
of the Poisson point process to eliminate the $\sup_{(a,t_0)\in\R^2}$
and put $a=t_0=0$.

Here, we firstly introduce an FKG inequality for probability measures on
the path space, which will play an important role in our proofs.
Let $\xi=\xi^{(0,0)}$ be the random path starting at the origin
defined in (\ref{13}); denote by $\bar\Pi$ the space of paths where
$\xi$ takes value. We define a partial order ``$\preceq$" on
$\bar\Pi$ as follows. Given $\pi_1,\pi_2\in \bar\Pi$,
\begin{equation}\label{1023}
\pi_1\preceq\pi_2 \,\,\mbox{ if and only if }\,\,\pi_1(t)-\pi_1(s)\le\pi_2(t)-\pi_2(s)
\,\,\mbox{ for all }\,\,t\ge s\ge 0.
\end{equation} 
Define increasing events in $\bar\Pi$ as usual. Denote
by $\mu_\xi$ the distribution of $\xi$ on $\bar\Pi$.

\begin{lemma}[{\bf FKG Inequality}]
  \label{2.10} 
  $\mu_\xi$ satisfies the FKG inequality, namely, for any increasing events
  $A,B\subseteq\bar\Pi$, $\mu_\xi(A\cap B)\ge \mu_\xi(A)\mu_\xi(B).$
\end{lemma}

\vskip 2mm \noindent{\bf Proof.} Let $Z_i,i\in\N$ be an i.i.d.~ family of
random variables with uniform distribution on the interval $[-r,r]$. 
Let $\{N(t):t\ge 0\}$
be a one-dimensional Poisson process with parameter $2r\la$. 
Assume that $\{N(t): t\ge0\}$ is independent of $Z_i,i\in\N$. 
Define the random process $\{Y(t):t\ge 0\}$ as 
$$Y(t)=\Dp\sum_{i=0}^{N(t)} Z_i,\ \ \forall \ t\ge 0,$$ where $Z_0\equiv 0$. 
Then $Y$ has the same distribution as $\xi$. 

Now, for any given $n\in \N$, we define a {\it discrete} time random walk 
$Y_n$ such that $Y_n$ converges in distribution to $Y$ as 
$n\rightarrow\infty$. Let $J_{n,i},i\in\N$ be an i.i.d.~family of
random variables such that $\P(J_{n,i}=1)=1-\P(J_{n,i}=0)=\frac {2r\la}n$; 
let $Z_{n,i}=Z_i\cdot J_{n,i}$. Let $l_n:=1/n$ be the unit length of time. 
Define $Y_n$ as 
$$Y_n(t)=\sum_{i=0}^{\lfloor nt \rfloor} Z_{n,i},\ \ \forall \ \ t\ge 0,$$ 
where $Z_{n,0}\equiv 0$ and $\lfloor nt \rfloor$ is the integer part of $nt$. 

For any given $t\ge 0$, define the random variable 
$$N_n(t):=|\{m: Z_{n,m}\not=0 \ {\rm and}\ 1\le m \le nt\}|.$$ 
Noticing that $N_n(t)$ converges in distribution to $N(t)$ as 
$n\rightarrow\infty$, it is straightforward 
to check that $Y_n$ converges in distribution to $Y$ as 
$n\rightarrow\infty$. Denote by $\mu_{_{Y_n}}$ 
the distribution of $Y_n$.

Considering the configuration space $\Omega_n:=[-r,r]^\N$, let
$\mu_n=\Pi_{_{i\in\N}}\mu_{_{Z_{n,i}}}$ be a product probability measure 
on $\Omega_n$, where
$\mu_{_{Z_{n,i}}}$ is the distributions of $Z_{n,i}$. Define the map
$\psi_n:\Omega_n\rightarrow\bar\Pi$ by
$\psi_n(z_1,\ldots,z_i,\ldots)=\pi$ such that
$$\pi(t)=\Dp\sum_{i=0}^{\lfloor nt \rfloor} z_i, \ \ \forall \ t\ge 0,$$ 
where $z_0\equiv 0$.

Obviously, $\psi_n$ is an increasing map from $\Omega_n$ to $\bar\Pi$ (under
the standard partial order on $\Omega_n$ and the order~(\ref{1023}) on $\bar\Pi$) and
$\psi_n(Z_{n,1},\ldots,Z_{n,i},\ldots)$ has the same distribution as $Y_n$.

For any increasing events $A,B$ contained in $\bar\Pi$, it is
straightforward to check that $\psi_n^{-1}(A)$, $\psi_n^{-1}(B)$ are
increasing and $\psi_n^{-1}(A\cap
B)=\psi_n^{-1}(A)\cap\psi_n^{-1}(B)$. So, \begin{equation}\label{k3}\begin{array}{ll} 
\mu_{_{Y_n}}(A\cap B)\!\!\!\!\!&=\mu_n(\psi_n^{-1}(A\cap B))
=\mu_n(\psi_n^{-1}(A)\cap\psi_n^{-1}(B))\\[2mm]
&\ge \mu_n (\psi_n^{-1}(A))\mu_n(\psi_n^{-1}(B))
=\mu_{_{Y_n}}(A)\mu_{_{Y_n}}(B). \end{array}\end{equation} Here we used the FKG inequality for 
the product measure $\mu_n$ on $\Omega_n$ in the third step. Thus we get
the FKG inequality for $\mu_{_{Y_n}}$. The Lemma follows immediately by 
taking $n \rightarrow \infty$.
$\square$

Let $\xi^{(0,0)}$, $\xi^{(\gamma, 0)}$, $\gamma\ge r$, be two random walks defined in (\ref{13}) from $S$. 
Note that $\{\xi^{(0,0)}$, $\xi^{(\gamma, 0)}\}$ makes of a simple system of coalescing random
walks starting at the space-time points $(0,0), (\gamma,0)$. Denote by $\Delta_\gamma$ the difference between 
the two walks, then, $\Delta_\gamma$ makes of a jump process in $[0,\infty)$ with absorbing state $0$:  
In $x\in [r, \infty)$, it has rates
$$(2r+x\wedge (2r))\lambda$$ and jump laws 
\begin{equation}\label{dis1}
\nu_x:=\frac{2r-x\wedge (2r)}{2r+x\wedge (2r)}\delta_{\{-x\}}+\frac{2(x\wedge (2r))}{2r+x\wedge (2r)}
U[r-x\wedge (2r),r],
\end{equation}
where $\delta_{\{-x\}}$ is the usual dirac measure and $U[r-x,r]$ is the uniform distribution on $[r-x,r]$. Note
that the process $\Delta_\gamma, \ \gamma\ge r$ will never visit $(0,r)$.

Denote by $\Xi_\gamma$, $\gamma\ge r$ the one-dimensional random walk starting at $\gamma$ 
with uniform rate $3r\lambda $ and uniform jumps in $[-r,r]$. Let 
$T=\inf \{t>0: \Xi_{\gamma}(t)\le 0\}$, 
firstly, we have

\begin{lemma}
\label{2.9} There exists a constant $c_1>0$ such that 
$\P(T>t)\le c_1/\sqrt{t}$ for any $t>0$, where $c_1$ 
depends on $r,\lambda$ and $\gamma$.
\end{lemma}

{\bf Proof.} This is certainly a well known result. On not having found a
reference, we briefly sketch an argument, which is rather standard. 
We give the result for the embedded chain of $\Xi_\gamma$, which we denote $\tilde\Xi_\gamma$.
The result for $\Xi_\gamma$ follows in a standard way. We use Skorohod's representation of
$\tilde\Xi_\gamma$ as
\begin{equation}\label{dec1}
\Xi_\gamma(n)=B(S_n),
\end{equation}
where $B(\cdot)$ is a Brownian motion of unit diffusion coefficient
started at $\gamma$, and $S_1,S_2,\ldots$
are the partial sums of i.i.d.~random variables with a positive exponential
moment (one of which has the distribution of the exit time of a standard
Brownian motion $B'$ from the interval $(-R,R)$, with $R$ a random variable
uniformly distributed in $(0,r)$ and independent of $B'$.
Let $\tilde T=\inf \{n>0: \tilde \Xi_{\gamma}(n)\le 0\}$.
Then $\tilde T>n$ implies that $M(S_n)>-r$, where $M(\cdot)$ is the running
maximum of $B(\cdot)$.
Thus
\begin{equation}\label{dec2}
\P(\tilde T>n)\leq \P(M(\mu n)>-r)+\P(S_n<\mu n),
\end{equation}
for all $\mu>0$, and the well known distribution of $M$,
along with the exponential moment condition of the increments of
$S_n$, via the standard large deviation estimate for sums of i.i.d.~random
variables with that condition, imply the claim of the lemma
for $\P(\tilde T>n)$. $\square$

\vspace{.3cm}

For the jump process $\Delta_{2r}$, 
let ${\cal T}=\inf \{t>0: \Delta_{2r}(t)=\xi^{(0,0)}(t)-\xi^{(2r, 0)}(t)=0\}$, we also have

\begin{lemma}
\label{2.11} There exists a constant $c_2>0$ such that $\P({\cal T}>t)\le
c_2/\sqrt{t}$ for any $t>0$, where $c_2$ depends on $r$ and $\lambda$. 
\end{lemma}

{\bf Proof.} Let $\Delta_\gamma^\prime$, $\gamma \ge r$, 
be the jump process in $[0,\infty)$ starting at $\gamma$: In $x\in [r, \infty)$, it has uniform
rate $3r\lambda$ ( be independent of $x$) and has jump laws given in (\ref{dis1}) (dependent on $x$). Note that
all  $\Delta_\gamma^\prime$, $\gamma \ge r$ can be coupled together such that, 
for any $\gamma\le\gamma'$, 
\begin{equation}\label{dis2}
\P(\Delta_{\gamma}^\prime(t)\le \Delta_{\gamma'}^\prime(t), \ \forall \ t\ge 0)=1. \end{equation}

On the other hand, let ${\cal T}_\gamma^\prime$ be the first hitting time to
$0$ of $\Delta_\gamma^\prime$, we have clearly
\begin{equation}\label{dis3}
\P({\cal T}>t)\le\P({\cal T}_\gamma^\prime >t), \ \forall \ t>0, \ \gamma\ge 2r.
\end{equation}
Now, we can prove the Lemma by two steps, in the first step, we couple $\Delta_r^\prime$ with a process 
$\bar\Delta$, which also has absorbing state $0$, 
such that ${\cal T}^\prime\le \bar{\cal T}$, the first hitting time to $0$ of $\bar\Delta$, with probability one. 
In the second step, we prove that $\P(\bar{\cal T}>t)\le c_2/\sqrt{t}$ for any $t>0$.

Step 1. For any $x\in [r,2r)$, let $p_x$ be the probability of the process $\Delta_x^\prime$ jump to $0$ in two steps, 
it is straightforward to calculate by (\ref{dis1}) that
$$p_x=\frac{2x}{2r+x}\int^{2r-x}_{r-x}\frac 1x \frac{2r-y}{2r+y}dy=\frac{4r}{2r+x}\ln \frac{4r-x}{3r-x}-\frac r{2r+x}.$$
Let $p=\inf_{r\le x<2r}p_x$, obviously, $p>0$, and let $q=1-p$.

Let us consider the behavior of the jump process $\Delta_{x}^\prime$ for some $x\in [r,2r)$ in 
the coming two clock ticks. Obviously, it should behave as the following three cases. 
a) jumps to $0$ at the first clock tick with probability $\frac {2r-x}{2r+x}$ and then stays there; 
b) stays in $[r,2r)$ at the first colck tick and jumps to $0$ at the second clock tick, 
note that the probability of this case is $p_x$ given above; 
c) does not jump to $0$ at both the two clock ticks. Note that in case c) the process should stay in $[r, 4r)$. 

Let $\bar \Delta$ be a jump process in $[0,\infty)$ with absorbing
state $0$ as follows. It starts at $4r$ and behaves as $\Delta_{4r}^{\prime}$ in $[2r,\infty)$. Once it hits $[r,2r)$, 
assume that it be in some $x\in[r,2r)$, it stays at $x$ and waits two clock ticks (with rate $3r\lambda$, 
the same rate as the process $\Delta_{4r}^{\prime}$) 
and then, either jumps to $0$ with probability $p$ or jumps back to $4r$ with probability $q$.

By the comments and setting in the last two paragraphs, 
one may couple $\Delta_{4r}^\prime$ and $\bar \Delta$ together
such that 
\begin{equation}\label{dis4}
\P({\cal T}_{4r}^\prime\le \bar{\cal T})=1.
\end{equation}
Note that the above coupling only guarantees (\ref{dis4}) but the stochastical domination between
$\Delta_{4r}^\prime$ and $\bar \Delta$.

Step 2. Let $T_{[r,2r)}$ be the first hitting time to $[r,2r)$ of $\Delta_{4r}^\prime$, let $\tau_1,\tau_2$, 
independent of $T_{[r,2r)}$, be two i.i.d. waiting times with rate $3r\lambda$. 
Let $\bar T=T_{[r,2r)}+\tau_1+\tau_2$ and $\bar T_i, i\in\N$ be a series of independent copies of $\bar T$. 
By the definition of the process $\bar \Delta$, $\bar {\cal T}$ has the same distribution as
\begin{equation}\label{dis5}
\tilde{\cal T}:=\sum^{N(p)}_{i=1} \bar T_i, 
\end{equation}
where $N(p)$ is geometric(p), 
the geometric distribution with parameter $p$, and
independent of $\bar T_i, i\in\N$. At this point, the claim of the lemma 
should be pretty clear, and could be argued in the following way.
Consider $G(s):=\sum_{n\geq1} \P(\tilde{\cal T}>n)\, s^n,\,|s|<1$.
{}From the relationship of $G$ and the moment 
generating function of $\tilde{\cal T}$, together with~(\ref{dis5}), 
we readily get the following.
\begin{equation}\label{dis6}
G(s)\leq\mbox{const}/\sqrt{1-s} 
\end{equation}
for $|s|<1$, from which there follows
\begin{equation}\label{dis7}
n\,\P(\tilde{\cal T}>n)\,e^{-1}\leq G(1-1/n).
\end{equation}
(\ref{dis7}) and~(\ref{dis6}) now yield the lemma.
 $\square$

\vspace{.3cm}

Now, we begin to verify conditions $B_1$ and $B_2$. By Lemma~\ref{2.5}, it is
straightforward to get that, for both Poisson 
trees $X_\delta$ and $Y_\delta$, 
\be
\label{f3}
\Dp
\limsup_{n \to \infty}
\P(\eta_{_{\delta_n}}(0,t;0,\epsilon)\ge 2)
=2\phi(\epsilon/\sqrt{2t})-1, 
\end{equation} 
where $\delta_n$ is any sequence of positive numbers converging to
$0$ as $n\to\infty$,
$\eta_\delta=\eta_{_{X_\delta}}$ or $\eta_{_{Y_\delta}}$, and
$\phi(x)$ is the standard normal distribution function. 
By (\ref{1000}) we also have 
\be
\label{1001}
\Dp
\limsup_{n \to \infty}
\P(\bar\eta_{_{\delta_n}}(0,t;0,\epsilon)\ge 2)
=2\phi(\epsilon/\sqrt{2t})-1. 
\end{equation} 
This gives $B_1$ and $B_1'$. 

Verifying $B_2'$ for the Poisson web $X_\delta$ is equivalent to checking that for any $t>0$
\begin{equation}\label{x1}
\epsilon^{-1}\limsup_{N \to \infty}
\P(\bar\eta_{_{X_1}}(0,tN;0,\epsilon\sqrt N)\ge 3) \to 0\hbox{ as } \epsilon\to 0+,
\end{equation}
where $X_1$ is defined in (\ref{10}).

For that, fix $t>0$. On the Poisson field $S$ with parameter $\lambda=\lambda_0=\sqrt 3/6$,
choose $r=r_0=\sqrt 3$, and then 
define $X_1$ as in (\ref{10}). We first condition the probability in~(\ref{x1})
on the set of points of intersection, in increasing order, 
of the paths $\xi^s$, $s\in S$, with
$[0,\epsilon \sqrt N]$, denoted $\{K_1,\ldots,K_J\}$, where $J,K_1,\ldots,K_J$
are random variables, with $J$ an integer which can equal $0$ (in which case
set of intersection points is empty by convention).
We note that by the definition of $\xi^s$, $s\in S$, no two distinct $K_i$'s 
can be at distance smaller than $r_0$.
For $\{x_1,\ldots,x_n\}\subset[0,\epsilon \sqrt N]$, let
 $\xi_j:=\xi^{(x_j, 0)}, 1\le j\le n$ as in (\ref{13}). 
Let $\eta'=\eta'(x_1,\ldots,x_n)=|\{\xi_j(tN): 1\le j\le n\}|$
(conventioned to be $0$ if $\{x_1,\ldots,x_n\}=\emptyset$. 
Clearly, $J,K_1,\ldots,K_J$ depend only on the points of $S$ below
time $0$. Thus, since 
$\eta'$ depends only on the points of $S$ above and at time $0$
for all $\{x_1,\ldots,x_n\}\subset[0,\epsilon\sqrt N]$, 
given $J=n,K_1=x_1,\ldots,K_J=x_n$ the probability in~(\ref{x1})
equals 
\begin{equation}\label{x100}
\P(\eta'\ge 3). 
\end{equation}
We derive below an upper bound for~(\ref{x100})
which is independent of $\{x_1,\ldots,x_n\}$ (see~(\ref{x2a}-\ref{x2})).
First, we enlarge, if necessary, the set $\{x_1,\ldots,x_n\}$
to make sure that $x_1=0$, $x_n=\epsilon\sqrt N$, and $r_0\leq x_j-x_{j-1}\leq2r_0$.
This also ensures that  $n\leq\epsilon\sqrt N/r_0+1$, and the enlargement can
only increase~(\ref{x100}).

Now for the bound. If $\eta'\ge 3$, then there should be some $1\le j\le n-1$ 
such that $\xi_{j-1}(tN)<\xi_j(tN)<\xi_n(tN)$. Hence,
\begin{equation}\label{x2a}
\ba{rl}
&\P(\eta'\ge 3)\le \Dp\sum^{n-1}_{j=2}\P(\xi_{j-1}(tN)<\xi_j(tN)<\xi_n(tN))\\[3mm]
&=\Dp\sum^{n-1}_{j=2}\int_{\bar \Pi_j}\P(\xi_{j-1}(tN)<\xi_j(tN)<\xi_n(tN)|\xi_j=\pi)\mu_{_{\xi_j}}(d\pi)\\[3mm]
&=\Dp\sum^{n-1}_{j=2}\int_{\bar \Pi_j}\P(\xi_{j-1}(tN)<\xi_j(tN)|\xi_j=\pi)
\P(\xi_j(tN)<\xi_n(tN)|\xi_j=\pi)\mu_{_{\xi_j}}(d\pi),
\ea
\end{equation}
where $\bar \Pi_j$ is the state space of $\xi_j$, and $\mu_{_{\xi_j}}$ its
distribution. 
In the latter equality,
we used the 
independence of $\xi_{j-1}(tN)<\xi_j(tN)$ 
and $\xi_j(tN)<\xi_n(tN)$ conditioned on $\xi_j=\pi$. 

We argue below that $\P(\xi_{j-1}(tN)<\xi_j(tN)|\xi_j=\pi)$ 
decreases in $\pi$ and $\P(\xi_j(tN)<\xi_n(tN)|\xi_j=\pi)$ increases in $\pi$.
This and the FKG Inequality for 
$\mu_{_{\xi_j}}$ (Lemma~\ref{2.10}) imply that the right hand side
of~(\ref{x2a}) 
is bounded above by
\begin{equation}\label{x2}
\ba{rl}
&\Dp\sum^{n-1}_{j=2}\int_{\bar
  \Pi_j}\P(\xi_{j-1}(tN)<\xi_j(tN)|\xi_j=\pi)\mu_{_{\xi_j}}(d\pi)
\\[3mm] &\hskip 10mm \cdot 
\Dp\int_{\bar \Pi_j}\P(\xi_j(tN)<\xi_n(tN)|\xi_j=\pi)\mu_{_{\xi_j}}(d\pi)\\[3mm]
&=\Dp\sum^{n-1}_{j=2}\P(\xi_{j-1}(tN)<\xi_{j}(tN))\P(\xi_j(tN)<\xi_n(tN))\\[3mm]
&\le \Dp\sum^{n-1}_{j=2}\P(\xi_{j-1}(tN)<\xi_j(tN))\P(\xi_0(tN)<\xi_{n}(tN))\\[3mm]
&\le(n-1)\P(\xi^{(0,0)}(tN)<\xi^{(2r_0,0)}(tN))
\P(\xi^{(0,0)}(tN)<\xi^{(\epsilon\sqrt N,0)}(tN))\\[3mm]
&\le\frac{\epsilon\sqrt N}{r_0}\,\P({\cal T}>tN)\,\P({\cal T}_{\epsilon, N} >tN),
\ea
\end{equation}
where ${\cal T}$ is the time when
$\xi^{(0,0)}$ and $\xi^{(2r_0,0)}$ meet and coalesce,
and ${\cal T}_{\epsilon, N}$ is the analogue time for $\xi^{(0,0)}$ and $\xi^{(\epsilon\sqrt N,0)}$.

To see that $\P(\xi_{j-1}(tN)<\xi_j(tN)|\xi_j=\pi)$ 
decreases in $\pi$, it is enough to consider $\pi_1\preceq\pi_2$ such that
$\pi_2(s)=\pi_1(s)$ for $0\leq s<t_0$ and $\pi_2(s)=\pi_1(s)+a$ for $s\geq t_0$,
for some $0<t_0<tN$ and $0<a<r_0$. As observed by Rongfeng Sun (see Acknowledgements
below), we can couple $\P(\,\cdot\,|\xi_j=\pi_1)$ and $\P(\,\cdot\,|\xi_j=\pi_2)$ by
shifting all the Poisson points of $S$ at and above time $t_0$ by $a$ units to
right, keeping the remaining points still. Since this operation preserves
the event $\{\xi_{j-1}(tN)<\xi_j(tN)\}$, we get the claimed monotonicity
of $\P(\xi_{j-1}(tN)<\xi_j(tN)|\xi_j=\pi)$. The argument for the monotonicity of
$\P(\xi_j(tN)<\xi_n(tN)|\xi_j=\pi)$ is similar.

Now, let us consider the items at the last line of equation (\ref{x2}). 
By Lemma 2.5, 
$$\limsup_{N \to \infty}\P({\cal T}_{\epsilon,N}>tN)=\P({\cal T}_{\epsilon, B}>t),$$ 
where ${\cal T}_{\epsilon, B}$ is the time when
two i.i.d.~Brownian motions starting at the same time at distance $\epsilon$
apart meet and coalesce. Thus the latter probability is an $O(\epsilon)$ for
every $t>0$ fixed. 
By Lemma~\ref{2.11}, $\P({\cal T}>tN))\le c_2/\sqrt{tN}$. These estimates imply that 
\begin{equation}\nonumber
\limsup_{N \to \infty} \P(\bar\eta_{_{X_1}}(0,tN;0,\epsilon\sqrt N)\ge 3)\le 
\limsup_{N \to \infty}
\frac{\epsilon\sqrt N}{r_0}\,\P({\cal T}>tN)\,\P({\cal T}_{\epsilon,N}>tN)
=O(\epsilon^2),
\end{equation}
and we get $B_2'$ for $X_\delta$. 


Finally, we verify $B_2'$ for Poisson web $Y_\delta$. Fix $t>0$. Given $\epsilon>0$ and 
$\{\delta_i>0, i\in\N\}$ such that 
$\delta_i\to 0$ as $i\to\infty$. 

Given $i$, on the Poisson field $S$ with parameter $\lambda(\delta_i)=\delta_i^{-1}$, 
choose $r=r(\delta_i)=(3\delta_i/2)^{1/3}$, and then 
define $Y_{\delta_i}=X( \lambda(\delta_i),r(\delta_i))$ as in (\ref{12}). As we did for $X_\delta$, 
conditioned on the set of points of intersection of the paths $\xi^s$, $s\in S$, with
$[0,\epsilon]$, define $\eta''$ be the number of the corresponding remaining paths in time $t$.
An analogous procedure gives that
$$\P(\eta''\ge 3)\le \frac{\epsilon}{r(\delta_i)}\,
\P({\cal T}_{r(\delta_i)}>t)\,\P({\cal T}_{\epsilon,i}>t),$$
where ${\cal T}_{r(\delta_i)}$ is the time when $\xi^{(0,0)}$ and $\xi^{(2r(\delta_i),0)}$ meet 
and coalease, and ${\cal T}_{\epsilon,i}$ is the analogue time for $\xi^{(0,0)}$ and $\xi^{(\epsilon,0)}$. 
By Lemma~\ref{2.5}, for any fixed $t$, 
$\limsup_{i \to \infty}\P({\cal T}_{\epsilon,i}>t)=O(\epsilon)$, so, to verify $B_2'$ for the Poisson web $Y_\delta$,
it is 
sufficient to prove that
\begin{equation}\label{y1}
\P({\cal T}_{r(\delta_i)}>t)=O(r(\delta_i)).
\end{equation}

As random walks in $[0,\infty)$, $\xi^{(0,0)}$ and $\xi^{(2r(\delta_i),0)}$ 
have rate $2r(\delta_i)\lambda(\delta_i)=2(3/2)^{1/3}\cdot \delta_i^{-2/3}$ 
(which tends to $\infty$ as $i\to\infty$). 
Let $\{\xi^{(0,0)'},\xi^{(r(\delta_i),0)'}\}$ 
be another (coalescing) random walk system, and the only difference from 
$\{\xi^{(0,0)},\xi^{(r(\delta_i),0)}\}$ is that, as single walks,
$\xi^{(0,0)'}$ and $\xi^{(2r(\delta_i),0)'}$ have unit rate. 
Let ${\cal T}^\prime_{r(\delta_i)}$ be the corresponding coalescence time 
for $\{\xi^{(0,0)'},\xi^{(r(\delta_i),0)'}\}$. 
It is clear that
\begin{equation}\label{y2}
\P({\cal T}_{r(\delta_i)}>t)=
\P({\cal T}^{\prime}_{r(\delta_i)}>t[2r(\delta_i)\lambda(\delta_i)])
=\P({\cal T}^{\prime}_{r(\delta_i)}>2(3/2)^{1/3}t\delta_i^{-2/3}).
\end{equation}
Using Lemma~\ref{2.11} for ${\cal T}^{\prime}_{r(\delta_i)}$, 
we get (\ref{y1}), and then get $B_2'$ for $Y_\delta$.

\section*{Acknowledgements}
This work was begun when one of us (X.-Y.~W.) was visiting the Statistics
Department of the Institute of Mathematics and Statistics of the
University of S\~ao Paulo. He is thankful to the probability group
of IME-USP for hospitality. We would like to thank
A.~Sarkar for introducing us to drainage networks, and C.~M.~Newman
for comments on the background paragraph of earlier versions.
We much thank Rongfeng Sun for pointing out important incorrections
on Subsection 2.2 of an earlier version, as well as for subsequent 
discussions on our fixing of them.

\end{document}